\documentclass[11pt]{article} 

\usepackage[latin1]{inputenc}
\usepackage[T1]{fontenc}
\usepackage{textcomp,graphicx,epsfig}
\usepackage{latexsym,amssymb,amsmath,amsthm}
\usepackage{fancyhdr,rotating,makeidx,eurosym}
\usepackage{lmodern,url,xspace,hyperref,layout}

\newcommand\RR{\mathbb{R}}

\newcommand\dd{\mathrm{d}}
\newcommand\ee{\mathrm{e}}

\newtheorem{definition}{Definition}[section]
\newtheorem{lemma}{Lemma}[section]
\newtheorem{theorem}{Theorem}[section]

\begin{document}

\title{Bounds for left and right window cutoffs\\
\textit{\small{Dedicated to the memory of B\'eatrice Lachaud}}} 


\author{%
Javiera Barrera\footnote{Facultad de Ingeniería y Ciencias, 
Universidad Adolfo Ib\'a\~{n}ez. 
Av. Diagonal las Torres 2640  Pe\~{n}alol\'en, Santiago,  Chile.
\texttt{javiera.barrera@uai.cl}}
\and 
Bernard Ycart\footnote{Laboratoire Jean Kuntzmann, Univ. Grenoble-Alpes,
51 rue des Math\'ematiques 38041 Grenoble cedex 9, France.
\texttt{Bernard.Ycart@imag.fr}}}

\maketitle

\abstract{The location and width of the time window
in which a sequence of processes converges to equilibrum 
are given under conditions of exponential convergence. The
location depends on the side: the
left-window and right-window cutoffs may have different locations.
Bounds on the distance to equilibrium are given for both
sides. Examples prove that the bounds are tight.
}



\section{Introduction}
The term ``cutoff'' was introduced by Aldous and Diaconis 
\cite{AldousDiaconisShuffle}, to describe the phenomenon of abrupt
convergence of shuffling Markov chains. Many families of
stochastic processes have since been shown to have similar
properties: see \cite[Chap.~8]{LevinPeresWilmer2006} for an introduction
to the subject, \cite{Saloff-CosteRW} for a review of
random walk models in which the phenomenon occurs, and 
\cite{ChenSaloff-Coste2} for an overview of the theory. Consider 
a sequence of stochastic processes in
continuous time, each converging to a stationary distribution. 
Denote by $d_n(t)$ the distance between the distribution at time $t$ of
the $n$-th process and its stationary distribution, the `distance'
having one of the usual definitions (total variation,
separation, Hellinger, relative entropy, $L^p$, etc.). 
The phenomenon can be expressed at three increasingly sharp levels 
(more precise definitions will be given in section \ref{sec:distandcutoff}).
\begin{enumerate}
\item
The sequence has a cutoff at $(t_n)$ if 
$d_n(ct_n)$ tends to the maximum $M$ of the distance if
$c<1$, to $0$ if $c>1$. 
\item 
The sequence has a window cutoff at $(t_n,w_n)$ if
$\liminf d_n(t_n+cw_n)$ tends to $M$ as $c$ tends to $-\infty$, and
$\limsup d_n(t_n+cw_n)$ tends to $0$ as $c$ tends to $+\infty$.
\item
The sequence has a profile cutoff at $(t_n,w_n)$ with profile $F$ if 
$F(c)=\lim d_n(t_n+cw_n)$ exists for all $c$, 
and $F$ tends to $M$ at $-\infty$, to $0$
at $+\infty$.
\end{enumerate}
There are
essentially two ways to interpret the cutoff time $t_n$: as a mixing
time \cite[Chap. 18]{LevinPeresWilmer2006}, or as a hitting time
\cite{MartinezYcart01}. For
samples of Markov chains, the latter interpretation can be used to
determine explicit online stopping times for MCMC algorithms
\cite{MonteCarloYcart,Lachaud05,LachaudYcart06,DiedhiouNgom09}.

Sequences of processes for which an explicit profile can be determined 
are scarce. The first example of a window cutoff concerned the random
walk on the hypercube for the total variation distance; 
it was treated by Diaconis and Shahshahani 
shortly after the introduction of the notion \cite{Diaconis87}. It was 
soon precised into a profile cutoff by Diaconis, Graham, and Morrison 
\cite{Diaconis90}. Cutoffs for random walks on more general
products or sums of 
graphs have been investigated in \cite{cutoffgrafosYcart}, and
more recently by Miller and Peres \cite{MillerPeres12}.
Random walks on the hypercube can be interpreted as 
samples of binary Markov chains. Diaconis et al.'s results were
generalized to samples of continuous and discrete time finite state
Markov chains for the chi-squared and total variation distance in
\cite{Ycart99}, then to samples of more general processes, for four
different distances in \cite[section 5]{cutoffjBbLbY} (see also
\cite[Chap. 20]{LevinPeresWilmer2006}).
Other examples of profile cutoffs include the riffle shuffle for the 
total variation distance \cite{BayerDiaconis}, and birth and
death chains for the separation distance
\cite{DiaconisSaloff-costeB-D} or the total variation distance
\cite{perescutoff-BD}.
When the maximum $M$ of the distance is $1$ (total variation,
separation), the profile $F$ decreases from $1$ to $0$. Thus it can be
seen as the survival function of some probability distribution on the
real line. A Gaussian distribution has been found for the  
riffle shuffle with the total variation distance 
\cite[Theorem~2]{BayerDiaconis} or for some birth and death chain with
the separation distance \cite[Theorem 6.1]{DiaconisSaloff-costeB-D}.
A Gumbel distribution has been found for samples of finite Markov
chains and the total variation distance
\cite{Diaconis90,Ycart99}. For the Hellinger, chi-squared, or relative entropy
distances, other profiles were obtained in \cite{cutoffjBbLbY}. 

Explicit profiles are usually out of
reach, in particular for the total variation distance: only a window
cutoff can be hoped for. 
However the definition above, which
is usually agreed upon (\cite[Definition 2.1]{ChenSaloff-Coste2} 
or \cite[p.~218]{LevinPeresWilmer2006}), may
not capture the variety of all possible situations. As will be shown
here, the location of a left-window cutoff should be distinguished
from that of a right-window cutoff: see Figure 18.2, p.~256 of 
\cite{LevinPeresWilmer2006}.  
The main result of this note,  Theorem \ref{th:cutoffgene},
expresses the characteristics of the left and right windows 
in terms of a decomposition into exponentials of
the distances $d_n(t)$. It refines some of the  results in 
Chen and Saloff-Coste \cite{ChenSaloff-Coste3}, in particular
Theorem 3.8. 
Explicit bounds on the distance to equilibrium are given. They 
are proved to be tight, using examples
of cutoffs for Ornstein-Uhlenbeck processes (see Lachaud
\cite{Lachaud05}).

The paper is organized as follows.
Section \ref{sec:distandcutoff} contains formal definitions and
statements. Examples are given in section
\ref{sec:wincutoffexamples}. Theorem \ref{th:cutoffgene}
is proved in section \ref{sec:proof}.   
\section{Definitions and statements}
\label{sec:distandcutoff}
For each positive integer $n$ a stochastic 
process $X_n=\{X_n(t)\,;\;t\geqslant 0\}$ is given. We assume that
$X_n(t)$ converges in distribution to $\nu_n$ 
as $t$ tends to infinity. The convergence is measured by
one of the usual distances (total variation, separation, Hellinger, 
relative entropy,  $L^p$, etc.),
the maximum of which is denoted by $M$ ($M=1$ for total variation and
separation, $M=+\infty$ for relative entropy, chi-squared\ldots). The distance
between the distribution of $X_n(t)$ and $\nu_n$ is 
denoted by $d_n(t)$.
\begin{definition}
\label{def:cutoff}
Denote by $(t_n)$ and $(w_n)$ two sequences of positive reals, 
such that $w_n=o(t_n)$. They will be referred to respectively as
\emph{location} and \emph{width}. The
sequence $(X_n)$ has:
\begin{enumerate}
\item
a left-window cutoff at $(t_n,w_n)$ if:
\[
\lim_{c\to -\infty} \liminf_{n\rightarrow\infty}
\inf_{t<t_n+cw_n}d_n(t) = M\;;
\]
\item
a right-window cutoff at $(t_n,w_n)$ if:
\[
\lim_{c\to +\infty} \limsup_{n\rightarrow\infty}
\sup_{t>t_n+cw_n}d_n(t) = 0\;;
\]
\item
a profile cutoff at $(t_n,w_n)$ with profile $F$ if:
\[
\forall c\in \RR\,,\; 
F(c) = \lim_{n\to\infty} d_n(t_n+cw_n)
\]
exists and satisfies:
\[
\forall c\in \RR\,,\; 0<F(c)<M
\quad\mbox{and}\quad
\lim_{c\to-\infty}F(c)=M\;,\quad
\lim_{c\to+\infty}F(c)=0\;.
\]
\end{enumerate} 
\end{definition}
If both left- and right-window cutoffs hold for the same location
$t_n$ and width $w_n$, then a $(t_n,w_n)$-cutoff holds in the sense of
Definition 2.1 in Chen and Saloff-Coste
\cite{ChenSaloff-Coste2}. The
location and width are not uniquely determined. Observe that if a
left-window cutoff holds at location $t_n$, it also holds at any
location $t'_n$ such that $t'_n\leqslant t_n$. Symmetrically, if a
right-window cutoff holds at location $t_n$, it also holds at any
location $t'_n$ such that $t'_n\geqslant t_n$. Moreover, if a cutoff
holds for width $w_n$, it also holds for any width $w'_n$ such
that $w'_n\geqslant w_n$. The location and width of a left-window
cutoff will be said to 
be optimal if for any $c<0$:
\[
\liminf_{n\rightarrow\infty}
\inf_{t<t_n+cw_n}d_n(t) < M\;.
\]
Those of a right-window cutoff are optimal if for any $c>0$:
\[
\limsup_{n\rightarrow\infty}
\sup_{t>t_n+cw_n}d_n(t) >0\;.
\]
This corresponds to strong optimality in the sense of 
\cite[Definition 2.2]{ChenSaloff-Coste2}. Of course, if a profile cutoff
holds, then the left- and right-window cutoffs hold at the same
location and width, which are optimal for both. Examples
will be given in section \ref{sec:wincutoffexamples}.   

Our main result relates the location and width of the left- and 
right-window cutoffs to the terms of a decomposition into exponentials of
the functions $d_n(t)$. From now on, we assume $M=+\infty$: the
distance is relative entropy, $L^p$ for $p>1$, etc. 
The result is expressed for a sequence of continuous
time processes, it could be written in discrete time, at the expense
of heavier notations.
\begin{theorem}
\label{th:cutoffgene}
Assume that for each $n$, there exist an increasing sequence of positive reals
$(\rho_{i,n})$, and a sequence of non negative reals $(a_{i,n})$ with $a_{1,n}>0$,
such that:
\begin{equation}\label{distance}
d_n(t)= \sum_{i=1}^{+\infty} a_{i,n}\, \ee^{-\rho_{i,n}t}\;.
\end{equation}
Denote by $A_{i,n}$ the cumulated sums of $(a_{i,n})$, truncated to
values no smaller than $1$. 
\[
A_{i,n}=\max \{1,a_{1,n}+\cdots+a_{i,n}\}\;.
\]
For each $n$, define:
\begin{equation} \label{def:tn}
t_n = \sup_{i} \frac{\log(A_{i,n})}{\rho_{i,n}}\;,
\end{equation}
\begin{equation} \label{def:wn}
w_n = \frac{1}{\rho_{1,n}}\;,
\end{equation}
\begin{equation}\label{def:rn}
r_n=w_n\left(\log(\rho_{1,n}t_n)-
\log(\log(\rho_{1,n}t_n))\right)\;. 
\end{equation}
Assume that:
\begin{enumerate}
\item  for $n$ large enough,
\begin{equation}\label{cond:tn}
0<t_n<+\infty\;,
\end{equation}
\item 
\begin{equation}\label{cond:rho1}
\lim_{n\rightarrow\infty} \rho_{1,n}t_n = +\infty\;,
\end{equation}
\item there exists a positive real $\alpha$ such that for $n$ large
  enough, and for all $i\geqslant 2$,
\begin{equation}\label{cond:alpha}
a_{i,n} \leqslant \alpha A_{i-1,n}\;.
\end{equation}
\end{enumerate}
Then $(X_n)$ has a left-window cutoff at $(t_n,w_n)$, a right-window
cutoff at $(t_n+r_n, w_n)$. More precisely:
\begin{equation}\label{teo:liminf}
\forall c<0\;,\quad \liminf_{n \to \infty}
d_n(t_n+c w_n) \geqslant \ee^{-c}\;,
\end{equation}
\begin{equation}\label{teo:limsup}
\forall c>0\;,\quad \limsup_{n \to \infty}
d_n(t_n+r_n+c w_n)\leqslant \ee^{-c}\;. 
\end{equation}
\end{theorem}
Conditions
(\ref{cond:tn}) and 
(\ref{cond:alpha}) are technical. Condition (\ref{cond:rho1}) is known
as Peres criterion: Chen and Saloff-Coste \cite{ChenSaloff-Coste2}
have proved that it implies cutoff for
$L^p$ distances with $p>1$, and given a counterexample for the 
$L^1$ distance. A consequence is that $w_n=o(t_n)$ as requested by
Definition \ref{def:cutoff}, and more precisely 
that $w_n=o(r_n)$ and $r_n=o(t_n)$.

A decomposition into exponentials of the distance to equilibrium such
as (\ref{distance}) holds for many processes: 
functions of finite state space Markov chains, functions of
exponentially ergodic Markov processes, etc. Assuming that the
decomposition only has non-negative terms is a stronger
requirement: see \cite[section 4]{ChenSaloff-Coste3}. It implies that
$d_n(t)$ is a decreasing function of $t$. We do not view it as a
limitation. Indeed, if (\ref{distance})
has negative terms, it can be decomposed as
$d_n(t)=d_n^+(t)-d_n^-(t)$, with:
$$
d_n^+(t) = \sum_{i=1}^{+\infty} \max\{a_{i,n},0\}\,\ee^{-\rho_{i,n} t}
\quad\mbox{and}\quad
d_n^-(t) = -\sum_{i=1}^{+\infty} \min\{a_{i,n},0\}\,\ee^{-\rho_{i,n} t}\;.
$$
Assume that Theorem \ref{th:cutoffgene} applies to both $d^+_n(t)$
and $d^-_n(t)$, leading to left-window cutoffs at $(t^+_n,w^+_n)$
and $(t^-_n,w^-_n)$, right-window cutoffs at $(t^+_n+r^+_n,w^+_n)$
and $(t^-_n+r^-_n,w^-_n)$. Since $d_n(t)$ is nonnegative, $t^-_n\leqslant
t^+_n$, $t^-_n+r^-_n\leqslant t^+_n+r^+_n$, and $w^-_n<w^+_n$.
The sequence $(X_n)$ has a
right-window cutoff, and (\ref{teo:limsup})
holds for $d_n$ with $(t_n+r_n,w_n)=(t_n^++r_n^+,w_n^+)$. 
Moreover, if $t_n^-+r_n^-=o(t^+_n)$ then the sequence $(X_n)$ has a
left-window cutoff, and 
(\ref{teo:liminf})
holds for $d_n$ with $(t_n,w_n)=(t_n^+,w_n^+)$.

Theorem 3.8 in
\cite{ChenSaloff-Coste3} contains a less tight assertion: it
describes a $(t_n,r_n)$-cutoff, which can be deduced from 
Theorem \ref{th:cutoffgene} above. However, it hides the fact that when
there is a (two-sided) window cutoff, the optimal width is no larger 
than $w_n$ thus strictly smaller than $r_n$. The latter quantity is
a correction bound on the location rather than a width: 
the optimal location may be anywhere between $t_n$ 
and $t_n+r_n$.
 
In the next section, sequences of processes having a profile cutoff at
$(t_n,w_n)$ or $(t_n+r_n,w_n)$, with profile $F(c)=\ee^{-c}$ will be
constructed, thus proving
that (\ref{teo:liminf}) and (\ref{teo:limsup}) are tight.

\section{Examples}
\label{sec:wincutoffexamples}
Several examples from the existing literature could be written
as particular cases of Theorem \ref{th:cutoffgene}: reversible Markov
chains for the $L^2$ distance  \cite{Ycart99,ChenSaloff-Coste3},
$n$-tuples of independent processes for the relative entropy distance
\cite{cutoffjBbLbY}, random walks on sums or products of graphs 
\cite{cutoffgrafosYcart}, samples of Ornstein-Uhlenbeck processes 
\cite{Lachaud05}. The objective of this section is not an extensive
review of possible applications, but rather the explicit
construction of some sequences illustrating the tightness of
(\ref{teo:liminf}) and (\ref{teo:limsup}), and the possible locations
of window cutoffs. We shall use here
the relative entropy distance, also called Kullback-Leibler divergence: 
if $\mu$ and $\nu$ are two probability measures
with densities $f$ and $g$ with respect to $\lambda$, then: 
\[
d(\mu,\nu)= \int_{S_\mu} f \log(f/g)\,\dd \lambda\;,
\]  
where $S_\mu$ denotes the support of $\mu$. The main advantage of
choosing that distance is its simplicity for dealing with tensor
products:
\[
d(\mu_1\otimes \mu_2,\nu_1\otimes\nu_2)= d(\mu_1,\nu_1)+d(\mu_2,\nu_2)\;.
\]
Let $a$ and $\rho$ be two positive
reals. Our building block will be a
one-dimensional  Ornstein-Uhlenbeck
process, denoted by $X_{a,\rho}$ (see Lachaud
\cite{Lachaud05} on cutoff for samples of
Ornstein-Uhlenbeck processes).  The process $X_{a,\rho}$ is a solution of
the equation:
\[
\dd X(t) = -\frac{\rho}{2} X(t) \,\dd t +\sqrt{\rho} \,\dd W(t)\;,
\]
where $W$ is the standard Brownian motion. 
The distribution of $X_{a,\rho}(0)$ is
normal with expectation $\sqrt{2a}$ and variance
$1$. It can be easily checked that the distribution of $X_{a,\rho}(t)$ is
normal with expectation $\sqrt{2a}\,\ee^{-\rho t/2}$ and variance
$1$. Therefore the (relative entropy) distance to equilibrium is:
\[
d(t)=a\,\ee^{-\rho t}\;.
\] 
Consider now two sequences $(a_n)$ and $(\rho_n)$ of
positive reals, and assume that $(a_n)$ tends to infinity. 
Theorem \ref{th:cutoffgene} applies to the sequence of processes 
$(X_{a_n,\rho_n})$ with $a_{1,n}=a_n$, 
$\rho_{1,n}=\rho_n$, and $a_{i,n}=0$ for $i>1$. The location and width are:
\[
t_n = \frac{\log(a_n)}{\rho_n}
\quad\mbox{and}\quad
w_n=\frac{1}{\rho_n}\;.
\]
The sequence has a profile cutoff at $(t_n,w_n)$ with profile
$F(c)=\ee^{-c}$. Indeed: 
\[
d_n(t_n+cw_n)= a_n\ee^{-(\rho_n t_n+c)} = \ee^{-c}\;.
\]
Hence
(\ref{teo:liminf}) is tight. 
For $\rho_n\equiv\rho$, $X_{a_n,\rho}$ is a Markov process with a fixed
semigroup, and an increasingly remote starting point: cutoff for such
sequences were studied in \cite{MartinezYcart01}. 

Using tuples of independent Ornstein-Uhlenbeck processes, one can
construct sequences $X_n$ for which the distance to equilibrium is any
finite sum of exponentials. Let $m_n$ be an integer. For
$i=1,\ldots,m_n$, let $a_{i,n}$ and
$\rho_{i,n}$ be two positive reals. Define the process $X_n$ as:
\[
X_n = \left(X_{a_{1,n},\rho_{1,n}},\ldots,X_{a_{m_n,n},\rho_{m_n,n}}\right)\;,
\]
where the coordinates are independent, each being an
Ornstein-Uhlenbeck process as defined above. The distance to
equilibrium of $X_n$ is:
\begin{equation}
\label{def:dnmn}
d_n(t) = \sum_{i=1}^{m_n} a_{i,n}\,\ee^{-\rho_{i,n} t}\;.
\end{equation}
Let $n$ be an integer larger than $1$.
Let $\beta_n$ be a real such that 
$0\leqslant \beta_n\leqslant 1$. 
Define:
\begin{equation}
\label{def:arho1}
a_{1,n} = \ee^{n}\;,\quad
\rho_{1,n}=
\frac{n}{1+\frac{\beta_n}{n}\log\left(\frac{n}{\log(n)}\right)}
\;,
\end{equation}
and for $i=2,\ldots,m_n=9^n$,
\begin{equation}
\label{def:arhoi}
a_{i,n} = \ee^{-n}\;,\quad
\rho_{i,n}=
\log(\ee^n+(i-1)\ee^{-n})
\;.
\end{equation}
The following notation is introduced for clarity:
\[
\ell_n = \log\left(\frac{n}{\log(n)}\right)\;.
\]
Using (\ref{def:tn}), (\ref{def:wn}), and (\ref{def:rn}), one gets:
\begin{equation}
\label{exp:tnwnrn}
t_n = 1+\frac{\ell_n\beta_n}{n}=\frac{n}{\rho_{1,n}}\;,\quad
w_n = \frac{t_n}{n}\;,\quad
r_n = \frac{t_n \ell_n}{n}=\ell_n w_n\;.
\end{equation}
\begin{lemma}
\label{lem:betan}
Let $d_n$ be defined by (\ref{def:dnmn}), with $a_{i,n}$ and
$\rho_{i,n}$ given by (\ref{def:arho1}) and (\ref{def:arhoi}).

Assume the following limit (possibly equal to $+\infty$) exists:
\begin{equation}
\label{limbetan}
\gamma=\lim_{n\to\infty} (1-\beta_n)\ell_n\;.
\end{equation}
Then:
\begin{equation}
\label{limdn}
\forall c\in \RR\;,\quad
\lim_{n\to\infty} d_n\left(t_n+(1-\beta_n)r_n+cw_n\right) = \ee^{-c}(1+\ee^{-\gamma})\;.
\end{equation}
\end{lemma}
A few particular cases are listed below. They 
illustrate the variety of possible
behaviors.
\begin{itemize}
\item 
$\beta\equiv 1$: a cutoff with
  profile $2\ee^{-c}$ occurs at $(t_n,w_n)$. 
\item 
$\beta_n\equiv \beta \in [0,1)$: a cutoff with
  profile $\ee^{-c}$ occurs at $(t_n+(1-\beta)r_n,w_n)$.
For $\beta= 0$, this proves that (\ref{teo:limsup}) is tight.
\item 
$\beta_n=(1+(-1)^n)/2$: a left-window cutoff occurs at $(t_n,w_n)$, a 
right-window cutoff
  at $(t_n+r_n,w_n)$. The locations and width are
  optimal.
\item 
$\beta_n=1-\gamma/\ell_n$, with $\gamma>0$: a cutoff
  with profile $\ee^{-c}(1+\ee^{\gamma})$ occurs at $(t_n,w_n)$.
\item 
$\beta_n=1-(2+(-1)^n)/\ell_n$: a $(t_n,w_n)$-cutoff occurs,
$t_n$ and $w_n$ are optimal. Yet no value of $c$ is such that
  $d_n(t_n+cw_n)$ converges: there is no profile.
\end{itemize}
\begin{proof}
The main step is the following limit.
\begin{equation}
\label{limdn1}
\lim_{n\to\infty} d_n\left(1+\frac{\ell_n}{n}+\frac{c}{n}\right) =
\ee^{-c}(1+\ee^{-\gamma})\;.
\end{equation}
In the sum defining $d_n$, let us isolate the first term:
$
d_n\left(1+\frac{\ell_n}{n}+\frac{c}{n}\right) =
D_1+ D_2\;,
$
with
\[
D_1=
a_{1,n}\exp\left(-\rho_{1,n}\left(
1+\frac{\ell_n}{n}+\frac{c}{n}\right)\right)
\;\mbox{and}\;
D_2 = 
\sum_{i=2}^{m_n}a_{i,n}\exp\left(-\rho_{i,n}\left(
1+\frac{\ell_n}{n}+\frac{c}{n}\right)\right)\;.
\]
The first term is:
$$
D_1=\exp\left(-\frac{(1-\beta_n)\ell_n+c}{t_n}\right)\;.
$$
Its limit is $\ee^{-(\gamma+c)}$ because $(1-\beta_n)\ell_n$ tends to
  $\gamma$ and $t_n$ tends to $1$. The second term is:
\[
D_2 = \sum_{i=2}^{+\infty}
\ee^{-n}\,\left(\ee^n+(i-1)\ee^{-n}\right)^{-\left(
1+\frac{\ell_n}{n}+\frac{c}{n}\right)}\;.
\]
Thus $D_2$ is a Riemann sum for the decreasing function 
$x\mapsto x^{-\left(1+\frac{\ell_n}{n}+\frac{c}{n}\right)}$.
Therefore,
\begin{equation}
\label{boundint}
\int_{\ee^n+\ee^{-n}}^{\ee^n+m_n\ee^{-n}}
x^{-\left(1+\frac{\ell_n}{n}+\frac{c}{n}\right)}\,\dd x 
<D_2<
\int_{\ee^n}^{\ee^n+(m_n-1)\ee^{-n}}
x^{-\left(1+\frac{\ell_n}{n}+\frac{c}{n}\right)}\,\dd x \;.
\end{equation}
Now:
\[
\frac{(\ee^n)^{-\left(\frac{\ell_n}{n}+\frac{c}{n}\right)}}
{\frac{\ell_n}{n}+\frac{c}{n}}
=\ee^{-c} \frac{\log(n)}{\ell_n+c}\;,
\]
which tends to $\ee^{-c}$. Moreover,
\[
\frac{(\ee^n+(m_n-1)\ee^{-n})^{-\left(\frac{\ell_n}{n}+\frac{c}{n}\right)}}
{\frac{\ell_n}{n}+\frac{c}{n}}
\leqslant \frac{n}{\ell_n+c} \left(\frac{m_n^{1/n}}{\ee}\right)^{-(\ell_n+c)}\;,
\]
which tends to $0$ for $m_n=9^n>\ee^{2n}$. So the upper bound in
(\ref{boundint}) tends to $\ee^{-c}$.
There remains to prove that the difference between the two integrals
tends to $0$. That difference is smaller than:
\[
\int_{\ee^n}^{\ee^n+\ee^{-n}}
x^{-\left(1+\frac{\ell_n}{n}+\frac{c}{n}\right)}\,\dd x
=
\left(
\frac{(\ee^n)^{-\left(\frac{\ell_n}{n}+\frac{c}{n}\right)}}
{\frac{\ell_n}{n}+\frac{c}{n}}
\right)
\left(1-
(1+\ee^{-2n})^{-\left(\frac{\ell_n}{n}+\frac{c}{n}\right)}
\right)\;.
\]
We have seen that the first factor tends to $\ee^{-c}$. The second
factor tends to $0$, hence the result.

Let us now deduce (\ref{limdn}) from (\ref{limdn1}).
Using (\ref{exp:tnwnrn}),
$$
1+\frac{\ell_n}{n}+\frac{c}{n}=t_n+(1-\beta_n)\frac{r_n}{t_n}+c\frac{w_n}{t_n}\;.
$$
Hence:
\begin{equation}
\label{limdn2}
\lim_{n\to\infty} d_n\left(t_n+(1-\beta_n)\frac{r_n}{t_n}+c\frac{w_n}{t_n}\right) =
\ee^{-c}(1+\ee^{-\gamma})\;.
\end{equation}
Let us write:
\[
t_n+(1-\beta_n)\frac{r_n}{t_n}+c\frac{w_n}{t_n}
=
t_n+(1-\beta_n)r_n+cw_n
-
((1- \beta_n)r_n+cw_n)
\left(\frac{\ell_n\beta_n}{nt_n}\right)\;.
\]
Therefore:
\begin{eqnarray*}
0&\leqslant&
d_n\left(t_n+(1-\beta_n)\frac{r_n}{t_n}+c\frac{w_n}{t_n}\right)
-d_n\left( t_n+(1-\beta_n)r_n+cw_n\right)\\
&\leqslant&
\left(\exp\left(\rho_{1,n}\left(((1-\beta_n)r_n+cw_n)
\frac{\ell_n\beta_n}{n}\right)-1\right)\right)
\,d_n\left( t_n+(1-\beta_n)r_n+cw_n\right)\\
&=&
\left(
\exp\left(
\frac{\ell_n^2(1-\beta_n)\beta_n+c\ell_n\beta_n}{nt_n}\right)-1\right)
\,d_n\left( t_n+(1-\beta_n)r_n+cw_n\right)\;.
\end{eqnarray*}
Hence the difference tends to $0$, 
since $\frac{\ell^2_n}{n}$ tends to $0$.
\end{proof}
\section{Proof of Theorem \ref{th:cutoffgene}}
\label{sec:proof}
Proofs of inequalities (\ref{teo:liminf}) and
(\ref{teo:limsup}) are given below.  
\begin{proof}[Proof of (\ref{teo:liminf})]$~$
Let $c$ be a negative real.
Fix $\epsilon$ such that $0<\epsilon<-c$. Using 
(\ref{def:tn}), define $i^*_n$ as: 
\begin{equation}
\label{def:istar}
i^*_n = \min\left\{\,i\,,\;
t_n-\epsilon w_n\leqslant 
\frac{\log(A_{i,n})}{\rho_{i,n}} \leqslant t_n\,\right\}\;. 
\end{equation}
From (\ref{cond:rho1}), $t_n+cw_n$ is positive for $n$ large enough. Then:
\begin{eqnarray*}
d_n(t_n+cw_n) &=& \sum_{i=1}^{+\infty}a_{i,n}\exp(-\rho_{i,n} (t_n+cw_n))\\
&\geqslant& \sum_{i=1}^{i^*_n}
a_{i,n}\exp(-\rho_{i,n}(t_n+cw_n))\\ 
&\geqslant&A_{i^*_n,n}\exp(-\rho_{i^*_n,n}(t_n+cw_n))\\[2ex]
&\geqslant&\exp((-\epsilon w_n-cw_n)\rho_{i^*_n})\\
&\geqslant&\exp((-\epsilon w_n-cw_n)\rho_{1,n})\\
&=&\ee^{-c-\epsilon}\;.
\end{eqnarray*}
 Since the inequality holds for all $\epsilon>0$, the result follows.
\end{proof}

\begin{proof}[Proof of (\ref{teo:limsup})]$~$
Let $c$ be a positive real. Our goal is to prove the following inequality.
\begin{equation}
\label{goal}
d_n(t_n+r_n+c w_n)\leqslant \ee^{-(r_n+cw_n)
  \rho_{1,n}}\frac{t_n}{r_n+cw_n}\left(\frac{r_n+cw_n}{t_n}+
  \ee^{C_n}\right)\;,
\end{equation}
where $C_n$ tends to $0$ as $n$ tends to infinity. Let us first check
that (\ref{goal}) implies (\ref{teo:limsup}). Observe that
$\frac{r_n+cw_n}{t_n}$ tends to $0$.
Using (\ref{def:wn}) and
(\ref{def:rn}):
\[
\ee^{-(r_n+cw_n)
\rho_{1,n}}\frac{t_n}{r_n+cw_n}=\ee^{-c}\frac{1}{1-\frac{\log(\log(t_n
\rho_{1,n}))+c}{\log(t_n \rho_{1,n})}}\; .
\]
By (\ref{cond:rho1}) the right-hand side tends to $\ee^{-c}$, hence
the result.

To prove (\ref{goal}), split the sum defining 
$d_n(t_n+r_n+c w_n)$ into two parts $S_1$ and $S_2$,
with:
\[
S_1= \sum_{i=1}^l a_{i,n}\exp(-\rho_{i,n}(t_n+r_n+cw_n))
\quad\mbox{and}\quad
S_2= \sum_{i=l+1}^{+\infty} a_{i,n}\exp(-\rho_{i,n}(t_n+r_n+cw_n))\;.
\]
Using the fact that the $\rho_{i,n}$ are increasing,
\begin{equation}
\label{S1}
S_1 \leqslant A_{l,n}\exp(-\rho_{1,n}(t_n+r_n+cw_n))\;.
\end{equation}
To bound $S_2$, the idea is the same as in the proof of (\ref{limdn}).
From (\ref{def:tn}), $\exp(-\rho_{i,n}t_n)\leqslant
A_{i,n}^{-1}$. Therefore:
\begin{equation}
\label{S2-1}
S_2 \leqslant \sum_{l+1}^{+\infty} a_{i,n}A_{i,n}^{-(1+(r_n+cw_n)/t_n)}\;. 
\end{equation}
The function $x\mapsto x^{-(1+(r_n+cw_n))/t_n}$ is decreasing, and its
integral from $l$ to $+\infty$ converges. The right-hand side of
(\ref{S2-1}) is a Riemann sum for that integral. Therefore:
\begin{equation}
\label{S2-2}
S_2 \leqslant \frac{t_n}{r_n+cw_n}
A_{l,n}^{-(r_n+cw_n)/t_n}\;.
\end{equation}
Consider first the particular case
$t_n=\frac{\log(A_{1,n})}{\rho_{1,n}}$, or equivalently
$A_{1,n}=\exp(t_n\rho_{1,n})$. Applying (\ref{S1}) and (\ref{S2-2})
for $l=1$ yields:
\begin{equation}
\label{eq:uppercutoffbound1}
d_n(t_n+r_n+c w_n) \leqslant
\ee^{-(r_n+cw_n)\rho_{1,n}}\frac{t_n}{r_n+cw_n}
\left(\frac{r_n+cw_n}{t_n}+1\right) \;,
\end{equation}
which is (\ref{goal}) for $C_n=0$. Otherwise,
$A_{1,n}<\exp(t_n\rho_{1,n})$. Let $\epsilon$ be such that 
$0<\epsilon<(t_n\rho_{1,n}-\log(A_{1,n}))/w_n$. The index $i^*_n$
defined by (\ref{def:istar}) is larger than $1$. The set of integers
$l$ such that $A_{l,n}<\ee^{\rho_{1,n}t_n}$, contains $1$ and is
bounded by $i^*_n$. Therefore, there exists
$l_n>1$ such that:
\begin{equation}
\label{def:ln}
A_{l_n-1,n} < \ee^{\rho_{1,n}t_n} \leqslant A_{l_n,n}\;.
\end{equation}
Applying (\ref{S1}) and (\ref{S2-2}) to $l=l_n-1$ yields:
\begin{eqnarray}\label{eq:uppercutoffbound2}
d_n(t_n+r_n+c w_n)  &\leqslant& \ee^{-(r_n+cw_n) \rho_{1,n}} +
\frac{t_n}{r_n+cw_n}
\exp\left(-\frac{r_n+cw_n}{t_n}\log A_{l_n-1,n}\right)\nonumber\\
&=& \ee^{-(r_n+cw_n) \rho_{1,n}} +
\frac{t_n}{r_n+cw_n}\exp\left(-(r_n+cw_n)\rho_{1,n}\frac{\log
    A_{l_n-1,n}}{\rho_{1,n}t_n}\right)\nonumber\\ 
&=&\ee^{-(r_n+cw_n)
  \rho_{1,n}}\frac{t_n}{r_n+cw_n}\left(\frac{r_n+cw_n}{t_n}+
  \ee^{C_n}\right)\;.
\end{eqnarray}
with
\begin{equation}
\label{def:Cn}
C_n=(r_n+cw_n)\rho_{1,n}\left(1-\frac{\log
      A_{l_n-1,n}}{\rho_{1,n}t_n}\right)\;.
\end{equation}
We must prove that $C_n$ tends to $0$.
By (\ref{def:wn}) and (\ref{def:rn}):
\begin{equation}\label{eq:r_n+cw_n}
 (r_n+cw_n)\rho_{1,n}=\log(\rho_{1,n}t_n)-\log\log(\rho_{1,n}t_n)+c\;.
\end{equation}
From (\ref{def:ln}):
\begin{equation}\label{eq:inlogA}
0< 1-\frac{\log (A_{l_n-1,n})}{\rho_{1,n}t_n}\leqslant
\frac{1}{\rho_{1,n}t_n}
\log\left(1+\frac{a_{l_n,n}}{A_{l_n-1,n}}\right)\;.
\end{equation}
Plugging (\ref{eq:r_n+cw_n}) and  (\ref{eq:inlogA}) 
into (\ref{def:Cn}), for $n$ large enough:
\[
0 < C_n \leqslant  \left(
\frac{\log(\rho_{1,n}t_n)-\log\log(\rho_{1,n}t_n)+c}{\rho_{1,n}t_n}\right)\,
\log\left(1+\frac{a_{l_n,n}}{A_{l_n-1,n}}\right)\;. 
\] 
By (\ref{cond:rho1}), the first factor of the
right-hand side tends to $0$. Moreover,
condition (\ref{cond:alpha}) entails that for $n$ large enough:
\[
 \log\left(1+\frac{a_{l_n,n}}{A_{l_n-1,n}}\right) < \log(1+\alpha) \;.
\]
Hence the result.
\end{proof}



\providecommand{\bysame}{\leavevmode\hbox to3em{\hrulefill}\thinspace}

\noindent
\textbf{Acknowledgements:} J. Barrera  was partially supported by
grants Anillo ACT88, Fondecyt  
  n\textsuperscript{o}1100618, and Basal project CMM (Universidad de Chile).
B. Ycart was supported by Laboratoire d'Excellence TOUCAN
  (Toulouse Cancer).
\end{document}